\documentclass[preprint,12pt]{elsarticle}
\newtheorem{thm}{Theorem}
\newtheorem{lem}[thm]{Lemma}
\newdefinition{rmk}{Remark}
\newproof{pf}{Proof} 
\newdefinition{dfn}{Definition}

\usepackage{amssymb}

\journal{Linear Algebra and its Applications}

\begin{document}

\begin{frontmatter}



\title{On the canonical structure of regular pencil  of singular matrix-functions}


\author{S.V.~Gaidomak}

\address{Institute for System Dynamics and Control Theory SB RAS, Lermontov  str., 134, Irkutsk, Russia 664033}

\begin{abstract}
The work is devoted to investigation of the canonical structure of regular, in domain of definition  ${\cal U}\subseteq {\mathbb R}^{m}$, the pencil of matrix-functions  $A(x)+\lambda B(x)$.  It is supposed that $\det A(x)\equiv 0$  and $\det B(x)\equiv 0\;\; \forall \ x \in {\cal U}$, and all roots of the characteristic equation $\det(A(x)+\lambda B(x))=0$ with $\lambda=\lambda (x)\in C^{p}({\cal U })$, are of constant multiplicity.

\end{abstract}


\end{frontmatter}
In the study of linear partial differential-algebraic equations of the form
$$
A(x)u_{t}(x)+B(x)\sum_{i=1}^{n}u_{x_{i}}(x)+C(x)u(x)=f(x),
$$
where coefficients $A(x)$, $B(x)$ and $C(x)$ are singular, at every point of the domain of definition ${\cal U}$, matrix of $C^{p}({\cal U})$ functions depending on many variables, it remains an important question about the possibility of bringing the pencil $A(x)+\lambda B(x)$ to  a canonical form with the help of suitable $p$ times differentiable linear transformation (see, for example~\cite{luch}-\cite{god}). In the papers~\cite{gai1}-\cite{gai3}, we considered the linear partial differential-algebraic equations with the pencil $A(x)+\lambda B(x)$, which has the structure of the ``rank-degree''. This is perhaps the most simple structure of the pencil $A(x)+\lambda B(x)$, where it is supposed that $\det A(x)\equiv 0$ and $\det B(x)\equiv 0\ \forall\ x\in {\cal U}$, with non-multiple finite and infinite elementary divisors of the pencil. The structure of the matrix-functions pencil known as ``rank-degree'' in the first time appears when studying of differential-algebraic systems~\cite{chist} and later it was extended by the author of~\cite{chist} to the case, when the pencil of matrix-functions depend on the several variables.

In this paper, we consider some canonical srtucture of regular pencil $A(x)+\lambda B(x)$, where it is supposed that $\det A(x)\equiv 0$ and $\det B(x)\equiv 0$ in ${\cal U}$, with finite and infinite elementary divisors, which has constant multiplicity equal no less than one in ${\cal U}$. As a result, we obtain sufficient conditions for the $p$-smooth equivalence of the matrix-functions pencil to this canonical form.

Throughout the paper, we use the following notations. Let $A(x)$ and $B(x)$  be an $n\times n$ matrices of ${C}^{p}({\cal U})$ real-valued functions of $m$ real variables $x\equiv (x^{1}, x^{2},\dots,x^{m})\in {\cal U}\subseteq {\mathbb R}^{m}$, where ${\cal U}$ is supposed to be a compact and simply connected domain in the ${\mathbb R}^{m}$. The expression of the form $A(x)+\lambda B(x)$, with $\lambda=\lambda(x)\in {C}^{p}({\cal U})$ is called  the matrix-functions pencil. One says, that the pencil under consideration is regular in ${\cal U}$, if exists a function $c=c(x)\in {C}^{p}({\cal U})$, for which the condition $\det(A(x)+c B(x))\neq 0\; \forall\ x\in {\cal U}$
holds~\cite{gant}. In what follows, we need in the following two definitions.

\begin{dfn}~\cite{verb0},~\cite{verb1},~\cite{verb2}
\label{de1}
Two square matrices of ${C}^{p}({\cal U})$ functions $A(x)$ and  $\tilde A(x)$ of order $n$ are called $p$-smoothly similar if there exists some matrix-function  $T(x)$ satisfying the following conditions:
\begin{enumerate}[(i)]
\item the elements of $T(x)$ belong to  ${C}^{p}({\cal U})$;
\item $T(x)$ is nonsingular in the domain of definition ${\cal U}$;
\item $T^{-1}(x)A(x)T(x)=\tilde A(x)$ in ${\cal U}$.
\end{enumerate}
Moreover, if  $T(x)$ is the unitary matrix-function, then one says that the matrix-functions $A(x)$ and  $\tilde A(x)$ are $p$-smoothly unitarily similar in ${\cal U}$~\cite{gant}.
\end{dfn}

\begin{dfn}
Two pencils of square ${C}^{p}({\cal U})$ matrix-functions of order $n$, say,  $A(x)+\lambda B(x)$ and $\tilde A(x)+\lambda \tilde B(x)$
is called $p$-smoothly equivalent, if there exists a pair of matrix-functions $P(x)$ and $Q(x)$,  independing on $\lambda$ and satisfying the following conditions:
\begin{enumerate}[(i)]
\item the elements of $P(x)$ and $Q(x)$ belong to  ${C}^{p}({\cal U})$;
\item  $P(x)$ and $Q(x)$ are nonsingular in ${\cal U}$;
\item  $P(x)(A(x)+\lambda B(x))Q(x)=\tilde A(x)+\lambda \tilde B(x)$ in ${\cal U}$.
\end{enumerate}
\end{dfn}

The investigation of the structure of pencils in the classical theory of matrices is commonly based on the properties of their similarity to the canonical forms. When  studying  the matrix-functions pencils, we keep this trend. 

Our paper is organized as follows. In the second section, we prove an auxiliary result concerning $p$-smooth similarity of matrix-functions with a single eigenvalue  identically equal to zero to some nilpotent matrix. In the third section, we prove a theorem on the $p$-smooth equivalence of pencil $A(x)+\lambda B(x)$ to some specified canonical form. This theorem is the main result of our paper. In conclusion of the paper, we give two examples of pencils to illustrate our theorem.

\section{Lemma on  the similarity of  matrix-functions}
\label{s1}

A well-known fundamental theorem of Shur and Toeplitz~\cite{lank} guarantees that every constant matrix $A$ is unitarily similar to a triangular matrix. Let us put the following question: whether this statement holds for matrix-functions? It turns out that for analytic matrix-functions of one variable the theorem of Shur and Toeplitz remains valid~\cite{gin}, but for analytic matrix-functions of several variables this theorem, in general, is not true. For example, consider the matrix-function with single eigenvalue identically equal to zero in the domain  of definition, namely,
$$
A(x)=\left (
\begin{array}{cc}
x_{1}(x_{1}+x_{2}) & -(x_{1}+x_{2})^{2}\\
x_{1}^2 & -x_{1}(x_{1}+x_{2})\\ 
\end{array}
\right ),
\ \ {\cal U}= {\mathbb R}^{2}.
$$ 
It is evident that all elements of $A(x)$ are analytic (polynomial) functions at every point of ${\cal U}$, but nevertheless the matrix-function $A(x)$ can not be cast to triangular form with the aid of a nonsingular analytic linear transformation, since the functions $\alpha (x)$ and $\beta(x)$ from lemma 1 of the paper~\cite{sil}, which must satisfy  the conditions:
$$
\alpha (x) x_{1}(x_{1}+x_{2})-\beta(x)x_{1}^2=0,\ \ \ \alpha^{2} (x)+\beta^{2}(x)=1
$$
are evidently not analytic at the point $(0,0)$. 

This example shows that $p$-smoothly unitary similar matrix-functions of several variables should satisfy additional conditions. In the case when the matrix-function have single eigenvalue being equal identically zero in ${\cal U}$,  additional condition to be required is the one of constant rank of the matrix-function in ${\cal U}$. Let us prove the following auxiliary statement.
\begin{lem}
Let $A(x)$ be an $n\times n$ matrix-function defined in ${\cal U}\subseteq {\mathbb R}^{m}$ satisfying the following conditions:
\begin{enumerate}[(i)]
\item the elements of $A(x)$  belong to the space  ${C}^{p}({\cal U})$;
\item  the $A(x)$ has  single eigenvalue being equal identically to zero in ${\cal U}$;
\item  $A(x)$ has constant rank in ${\cal U}$.
Then $A(x)$ is $p$-smoothly unitarily similar to some  nilpotent matrix-function ${\cal N}(x)$.
\end{enumerate}
\label{l1}
\end{lem}
\begin{pf}
The method of proof, in general, follows the line suggested in the papers \cite{lank} and \cite{gin}. Thus, there is no need to write down in detail all the proof of the lemma. Let us consider only  some important aspects of the proof. 

Let $X(x)$ be a right eigenvector of $A(x)$,  corresponding to the eigenvalue $\lambda \equiv 0$ in ${\cal U}$. Let us prove that the elements of eigenvector $X(x)$ belong to the  space ${C}^{p}({\cal U})$. We write down the equation for the right eigenvector
\begin{equation}
A(x)X(x)=0.
\label {eq2.1}
\end{equation}  
Since, by assumption, the rank of $A(x)$ does not depend on  $x$ in ${\cal U}$, then (see lemma~2, \cite{verb2}) there exists  a pair of ${C}^{p}({\cal U})$ matrix-functions $P(x)$ and $Q(x)$ nonsingular in ${\cal U}$, such that 
$$
P(x)A(x)Q(x)={\rm diag} \{E_{r},{\cal O}\},\ \ \mbox{where} \ \ r\equiv{\rm rank}(A(x)).
$$  
It is obvious that $r<n$. Let $\{e_i,\ i=\overline {1,n}\}$
be standard orthonormal basis of  linear space ${\mathbb R}^{n}$. It is easily seen that vector-function $X(x)=Q(x)e_{n}$ is the solution of equation~(\ref{eq2.1}).  Furthermore, the all components of vector $X(x)$ belong to the space ${C}^{p}({\cal U})$. Let us show now that it is possible to construct, attached to $X(x)$, the orthonormal system of linearly independent vectors
\begin{equation}
\hat X_{1}(x),\ \hat X_{2}(x),\ \dots,\ \hat X_{n}(x),
\label {eq2.2}
\end{equation}  
satisfying the following conditions:\footnote{Here $\|\cdot\|$ denotes euclidean norm.}
\begin{equation}
\hat X_{1}(x)=X(x)/\|X(x)\|,
\label {eq2.3}
\end{equation}   
\begin{equation}
\hat X_{i}(x)^{\top}\cdot\hat X_{j}(x)=0,\ \ \forall \ i\neq j,\ \ i,j=\overline{1,n}\ \ \mbox{è}\ \ \forall\ x\in {\cal U}, 
\label {eq2.4}
\end{equation}  
\begin{equation}
X_{i}(x)\in {C}^{p}({\cal U}),\ \ \|\hat X_{i}(x)\|\equiv 1, \ \ \forall\ i=\overline{1,n}.
\label {eq2.5}
\end{equation}  
Consider following set of vector-functions~\cite{gin}
\begin{equation} 
Z_{1}(x)=Q(x)e_{n},\ Z_{2}(x)=Q(x)e_{n-1},\ \dots,\ \ Z_{n}(x)=Q(x)e_{1}.
\label {eq2.6}
\end{equation}  
Since the matrix-function $Q(x)$ is nonsingular in ${\cal U}$, then~(\ref{eq2.6}) gives the set of linearly independent vector-functions ${\cal U}$. 
Applying the Gram-Schmidt procedure to~(\ref{eq2.6}), we obtain the orthogonal set of vectors
$$
X_{1}(x)=Z_{1}(x),
$$ 
$$
X_{2}(x)=Z_{2}(x)-\frac{X_{1}^{\top}(x)\cdot Z_{2}(x)}{X_{1}^{\top}(x)\cdot X_{1}(x)}X_{1}(x),
$$
\begin{equation}  
X_{3}(x)=Z_{3}(x)-\frac{X_{1}^{\top}(x)\cdot Z_{3}(x)}{X_{1}^{\top}(x)\cdot X_{1}(x)}X_{1}(x)-\frac{X_{2}^{\top}(x)\cdot Z_{3}(x)}{X_{2}^{\top}(x)\cdot X_{2}(x)}X_{2}(x),
\label {eq2.7}
\end{equation}  
$$
\dots\ \ \ \ \ \ \ \dots \ \ \ \ \ \ \dots \ \ \ \  \ \  \dots \ \ \ \ \ \ \dots\ \ \  \ \ \ \dots
$$
$$  
X_{n}(x)=Z_{n}(x)-\frac{X_{1}^{\top}(x)\cdot Z_{n}(x)}{X_{1}^{\top}(x)\cdot X_{1}(x)}X_{1}(x)-\cdots-\frac{X_{n-1}^{\top}(x)\cdot Z_{n}(x)}{X_{n-1}^{\top}(x)\cdot X_{n-1}(x)}X_{n-1}(x).
$$
By construction, each of vector-functions $X_{i}(x)$ is not equal to zero vector in ${\cal U}$. This means that
$$
X_{i}^{\top}(x)\cdot X_{i}(x)\neq 0,\ \ \ i=\overline{1,n},\ \ \forall\ x\in {\cal U}
$$ 
and the components of vector-functions~(\ref{eq2.7}) belong to the space ${C}^{p}({\cal U})$. Moreover, the vector-functions
$
\hat X_{i}(x)=X_{i}(x)/\|X_{i}(x)\|
$
satisfy required condition~(\ref{eq2.3})-(\ref{eq2.5}). Then, clearly, the matrix-function
$$
U(x)=(\hat X_{1}(x), \hat X_{2}(x),\dots,\hat X_{n}(x))
$$
is nonsingular in ${\cal U}$. Moreover, it satisfies the condition $U^{\top}(x)=U^{-1}(x)$, that is, $U(x)$ is the unitary matrix-functions.
It remains to prove that
$$
U(x)^{-1}A(x)U(x)={\cal N}(x),
$$ 
where ${\cal N}(x)$ is some nilpotent matrix-function. This part of the proof goes by induction and completely coincides with the proof of theorem Shur and Toeplitz~\cite{lank}. Therefore, we will not dwell on this part of the proof. So, the lemma is proved.
\end{pf}
It is worth to remark the third condition of the lemma is sufficient, but not necessary. We can give an example, when the rank of $A(x)$  is variable in the domain of definition ${\cal U}$, but nevertheless  $A(x)$ is $p$-smoothly similar to some nilpotent matrix-function.  
For example, consider the matrix-function
$$
A(x)=\left (
\begin{array}{ccc}
x_{1} & -1 & x_{1}(x_{1}+x_{2}) \\
x_{1}^{2} & -x_{1} & -x_{1}-x_{2}\\ 
0 & 0 & 0
\end{array}
\right )
,\ \ {\cal U}={\mathbb R}^{2}.
$$  
On the line given by the equation $x_{2}=-x_{1}$, the rank of $A(x)$ is equal to one, while outside this line the rank of $A(x)$ is two. Remark, that in  this example, there are not isolated points of change of the rank. In this case, we are able construct the unitary matrix-function $U(x)$ in ${\cal U}$. It is given by
$$
U(x)=\frac{1}{\sqrt{1+x_{1}^{2}}}
\left (
\begin{array}{ccc}
1 & x_{1} & 0\\
x_{1} & -1 & 0\\
0 & 0 & 1
\end{array}
\right )
$$
while the nilpotent matrix-function ${\cal N}(x)$ takes the following form:
$$
{\cal N}(x)=(1+x_{1}^{2})
\left (
\begin{array}{ccc}
0 & 1 & 0\\ 
0 & 0 & x_{1}+x_{2}\\ 
0 & 0 & 0
\end{array}
\right ).
$$

\section{The theorem on canonical structure of regular matrix-functions pencil}
\label{s3}
Let us prove the theorem on canonical form of regular matrix-functions pencil.
\begin{thm}
Let the following conditions be satisfied:
\begin{enumerate}[(i)]
\item the all roots of characteristic polynomial  $\det (A(x)+\lambda B(x))$  are real and of constant multiplicity in the domain of definition ${\cal U}$;
\item  the leading coefficient of polynomial  $\det (A(x)+\lambda B(x))$ is not identically equal to zero on ${\cal U}$;
\item  the ranks of $A(x)$ and $B(x)$  are independent of $x\in{\cal U}$ and less than dimension $n$.
\end{enumerate}
Then the pencil $A(x)+\lambda B(x)$ is $p$-smoothly equivalent to the  following canonical form:
\begin{equation}
 {\rm diag}\{E_{d}, {\cal M}(x),  E_{\hat l}\}+\lambda\ {\rm diag}\{{\cal J}(x), E_{l},  {\cal N}(x)\},
\label{eq3.1}
\end{equation}
where $E_{d}$ denotes identity matrix of order $d$; ${\cal M}(x)$ and  ${\cal N}(x)$ are some nilpotent matrices of  orders $l$ and $\hat l$, respectively, $\hat l=n-d-l$;
${\cal J}(x)={\rm diag}\{{\cal J}_{1}, {\cal J}_{2},\dots, {\cal J}_{k}\}$, where ${\cal J}_{i}$ for $i=\overline{1,k}$ are nonsingular  matrix-functions of orders $p_{i}$, respectively; $\sum_{i=1}^{k}p_{i}=d$;  every block ${\cal J}_{i}$ has single eigenvalue $-1/\lambda_{i}(x)$ in ${\cal U}$; $\lambda_{i}(x)$ for $i=\overline{1,k}$ being eigenvalues of the characteristic polynomial $\det (A(x)+\lambda B(x))$ are not equal to zero in ${\cal U}$.
\label{t1}
\end{thm}
\begin{pf}
It is obvious, that the coefficients of characteristic polynomial $\Delta (\lambda, x)=\det(A(x)+\lambda B(x))$ belong to the space ${C}^{p}({\cal U})$. 
Moreover, it is known, that the roots  $\lambda_{i}(x)$ for $i=\overline{1,k+1}$ belong to the space ${C}^{p}({\cal U})$ (see lemma~1,~\cite{verb1}). First condition of the theorem, require that the roots of $\Delta (\lambda, x)$ do not coincide in  ${\cal U}$, that is, $\lambda_{i}(x)\neq \lambda_{j}(x)\ \ \forall \ i\neq j$ and $\forall \ x\in {\cal U}$. The second condition of the theorem excludes the case, when $\Delta (\lambda, x)\equiv 0$ in ${\cal U}$, that is, the case of singular pencil. In these circumstances, there exists some function $c=c(x)\in {C}^{p}({\cal U})$, such that $c(x)\neq 0$ $\forall \ x\in {\cal U}$ for which  the condition
$$
\det(A(x)+c B(x))\neq 0\ \ \forall \ x\in {\cal U}
$$  
holds. The function $c(x)$ can be constructed, for instance, as the arithmetic mean of the two neighboring roots of characteristic polynomial $\Delta (\lambda, x)$.

Let  $A_{1}(x)\equiv A(x)+c B(x)$, then 
\begin{equation}
A(x)+\lambda B(x)=A_{1}(x)+(\lambda-c)B(x).
\label{eq3.2}
\end{equation}
Multiplying the equation~(\ref{eq3.2}) by the matrix $A_{1}^{-1}(x)$ on the left, yields
\begin{equation} 
A_{1}^{-1}(x)\left [ A(x)+\lambda  B(x)\right ]=E_{n}+(\lambda-c)A_{1}^{-1}(x)B(x).
\label{eq3.3}
\end{equation}
Let us consider for matrix-function $A_{1}^{-1}(x)B(x)$ its characteristic polynomial
$$
\tilde \Delta(\xi,x)=\det (A_{1}^{-1}(x)B(x)-\xi E_{n}),
$$
with $\xi=\xi(x)$ being some unknown function. Making use of elemetary properties of the determinant, we can write
\begin{equation}
\tilde \Delta(\xi,x)=(-1)^{n}\det A_{1}^{-1}(x) \det (\xi A(x)+(\xi c-1)B(x)).
\label{eq3.4}
\end{equation} 
Let us write down the characteristic polynomial $\Delta(\lambda,x)$ in the form
\begin{equation}
\Delta(\lambda,x)=\sum_{i=0}^{n}S_{i}(x)\lambda^{i},
\label{eq3.5}
\end{equation} 
where  the coefficients $S_{i}(x)$ are the sums of the all minors of order $n$,  composed of the $n-i$ rows of $A(x)$ and  the $i$ rows of $B(x)$. 
Clearly, $S_{0}(x)$ and $S_{n}(x)$ are determinants of $A(x)$ and $B(x)$, respectively.

In virtue of the third condition the theorem, the ranks of $A(x)$ and $B(x)$ are less than dimension $n$. Hence, $S_{0}(x)\equiv 0$ and $S_{n}(x)\equiv 0$ in  ${\cal U}$. Then, among the roots of the polynomial $\Delta(\lambda,x)$ we can always find the root which is identically equal to zero in  ${\cal U}$. It can be assumed, without loss of generality, that, say, $\lambda_{k+1}(x)$ be a zero root of multiplicity $l$, where $l\geq 1$. From what we said, it follows that (\ref{eq3.5}) takes the form
\begin{equation}
\Delta(\lambda,x)=\lambda^{l}\sum_{i=l}^{l+d}S_{i}(x)\lambda^{i-l},\ \ l+d\leq n,
\label{eq3.6}
\end{equation} 
where $l$ and $d$ is independent of $x$ in ${\cal U}$. Since the multiplicities of the roots of polynomial  $\Delta(\lambda,x)$ are constant in ${\cal U}$, then coefficients $S_{l}(x)$ and $S_{l+d}(x)$  are not equal to zero in ${\cal U}$.
In virtue of conditions the theorem,  (\ref{eq3.6}) is also specified as
\begin{equation}
\Delta(\lambda,x)=S_{l+d}(x)\lambda^{l} \prod_{i=1}^{k}(\lambda-\lambda_{i}(x))^{p_{i}},\ \mbox{where} \ \sum_{i=1}^{k}p_{i}=d.
\label{eq3.7}
\end{equation}  
Together with $\Delta(\lambda,x)$, we consider the polynomial 
$$
\hat \Delta (\mu,\lambda, x)=\det(\mu A(x)+\lambda B(x))
$$  
with some $\mu$ being some unknown function $\mu=\mu(x)$. Let us write down the latter as 
$$
\hat\Delta(\lambda,x)=\sum_{i=0}^{n}S_{i}(x)\lambda^{i}\mu^{n-i}.
$$
Making use of ~(\ref{eq3.6}) and~(\ref{eq3.7}) we get
\begin{equation}
\tilde\Delta(\lambda,x)=(-1)^{n}\epsilon (x)\xi^{\hat l}\left (\xi -\frac{1}{c}\right )^{l} \prod_{i=1}^{k}\left (\xi -\frac{1}{c-\lambda_{i}(x)}\right )^{p_{i}},
\label{eq3.9}
\end{equation}  
with
$$
\epsilon(x)=\det A_{1}^{-1}(x)S_{l+d}(x) c^{l}\prod_{i=1}^{k}(c-\lambda_{i}(x))^{p_{i}}.
$$
It is worth to remember that in the relation ~(\ref{eq3.9}) $c\neq \lambda_{i}(x)$ for all $x\in{\cal U}$ and for all $i=\overline{1,n}$. On the other hand, we can write
\begin{equation}
\tilde\Delta(\lambda,x)=(-1)^{n}\xi^{n}+O(\xi^{n-1}),
\label{eq3.10}
\end{equation}  
where $O(\xi^{n-1})$ stands for the sum of the terms containing degrees of $\xi$ less, than $n$.
Comparing relations~(\ref{eq3.9}) and~(\ref{eq3.10}), we conclude that $\epsilon(x)\equiv 1$. Thus, 
$$
\tilde\Delta(\lambda,x)=(-1)^{n}\xi^{\hat l}\left (\xi -\frac{1}{c}\right )^{l} \prod_{i=1}^{k}\left (\xi -\xi_{i}(x)\right )^{p_{i}}\ \mbox{with} \ \ \xi_{i}(x)=\frac{1}{c-\lambda_{i}(x)}.
$$
Since $\lambda_{i}(x)\neq \lambda_{j}(x)$ for all $i,j=\overline{1,n}$ and $c\neq 0$ in ${\cal U}$, then, according to the theorem from~\cite{sib} there exists the matrix-function $T(x)$, which satisfies all the conditions of the definition~\ref{de1} and the following relation:
$$
T^{-1}(x)A_{1}^{-1}(x)B(x)T(x)=\mbox{diag}\{ J_{1}(x),J_{2}(x),\dots, J_{k}(x), M(x), N(x)\},
$$ 
where $J_{i}(x)$ for $i=\overline{1,k}$ are nonsingular in ${\cal U}$ and  square matrix-functions  with pair-wise unequal eigenvalues $\xi_{i}$ for $i=\overline{1,n}$. According to the theorem from~\cite{sib}, $M(x)$  is the nonsingular in ${\cal U}$  matrix-function of order $l$ with  one eigenvalue equal to $1/c$, while $N(x)$ is the matrix-function of order $\hat l$ with  one eigenvalue being identically equal to zero in ${\cal U}$.

Multiplying the pencil from~(\ref{eq3.3}) by the $T^{-1}(x)$ and $T(x)$ on the left and on the right, respectively, we obtain
$$
E_{n}+(\lambda-c)\mbox{diag}\{J_{1}(x),J_{2}(x),\dots, J_{k}(x), M(x), N(x)\}
$$
$$
=\mbox{diag}\{E_{p_{1}}-cJ_{1}(x),E_{p_{2}}-cJ_{2}(x),\dots,E_{p_{k}}-cJ_{k}(x), E_{l}-c M(x), E_{\hat l}-c N(x)\}
$$
\begin{equation}
+\lambda\ \mbox{diag}\{J_{1}(x),J_{2}(x),\dots, J_{k}(x), M(x), N(x)\}.
\label{eq3.11}
\end{equation}
Since, the ranks of matrix-functions $B(x)$, $J(x)$ and $M(x)$ are constants in ${\cal U}$, then, in virtue of the known property of rank~\cite{gant}, the rank $N(x)$ is independent of $x$ in ${\cal U}$. Lemma~\ref{l1} says, that $N(x)$ is $p$-smoothly unitary similar to some nilpotent matrix-function $\hat  N (x)$. This means that there exists the nonsingular unitary matrix-functions $U(x)$,  satisfying the relation
$$
U^{-1}(x)N(x)U(x)= \hat N(x).
$$ 
Take the following matrix-function:
$$
\tilde U(x)=\mbox{diag}\{E_{d},E_{l},U(x)\}.
$$
Multiplying the pencil~(\ref{eq3.11}) by the $\tilde U^{-1}(x)$ and $\tilde U(x)$ on the left and on the right, respectively, we bring this pencil to the form
$$
\mbox{diag}\{E_{p_{1}}-cJ_{1}(x),E_{p_{2}}-cJ_{2}(x),\dots,E_{p_{k}}-cJ_{k}(x), E_{l}- c M(x), E_{\hat l}- c \hat N (x)\}
$$
\begin{equation}
+\lambda\ \mbox{diag}\{J_{1}(x),J_{2}(x),\dots, J_{k}(x), M(x), \hat N (x)\}.
\label{eq3.12}
\end{equation} 
Let us consider now the  following blocks $E_{p_{i}}-cJ_{i}(x)$ for $i=\overline{1,k}$ and the characteristic polynomials for each of them:  
\begin{equation}
\det(E_{p_{i}}-cJ_{i}(x)-\nu_{i}E_{p_{i}})=(-1)^{p_{i}}c^{p^{i}}\det \left (J_{i}(x)-\frac{1-\nu_{i}}{c} E_{p_{i}}\right)
\label{eq3.121}
\end{equation} 
By assumption, the function $c$ must not coincide with any root of characteristic polynomial $\Delta(\lambda,x)$ in ${\cal U}$. In particular, this means that $c\neq 0$ in ${\cal U}$. Since, the roots of (\ref{eq3.121}), namely, $\nu_{i}={-\lambda_{i}(x)}/{(c-\lambda_{i}(x))}$  are not equal to zero in $ {\cal U} $, then the matrix-functions $E_{p_{i}}-cJ_{i}(x)$ are  nonsingular in  ${\cal U}$. Multiplying the pencil~(\ref{eq3.12}) by the matrix-function
$$
\bar J(x)=\mbox{diag}\{\hat J(x), E_{l}, (E_{\hat l}-c\hat N(x))^{-1}\},
$$   
on the left whith
$$
\hat J(x)=\mbox{diag}\{(E_{p_{1}}-cJ_{1}(x))^{-1},(E_{p_{2}}-cJ_{2}(x))^{-1},\dots,(E_{p_{k}}-cJ_{k}(x))^{-1} \},
$$
we obtain the pencil
\begin{equation}
\mbox{diag}\{E_{d}, E_{l}- c M(x),  E_{\hat l}\}+\lambda\ \mbox{diag}\{{\cal J}(x), M(x),  {\cal N}(x)\},
\label{eq3.13}
\end{equation} 
where
$$
{\cal J}(x)=\mbox{diag} \{{\cal J}_{1}(x),{\cal J}_{2}(x),\dots ,{\cal J}_{k}(x)\}  
$$
with blocks ${\cal J}_{i}(x)=(E_{p_{i}}-cJ_{i}(x))^{-1}J_{i}(x)$  and 
$
{\cal N}(x)=(E_{\hat l}-c\hat N(x))^{-1}\hat N(x).
$

The matrix-function ${\cal N}(x)$ is nilpotent, since it is constructed as a product of triangular and nilpotent matrix-functions. Each  block ${\cal J}_{i}(x)$ in~(\ref{eq3.13})  is nonsingular in ${\cal U}$ and has unique eigenvalue $-1/\lambda_{i}(x)\ \forall \ i$. Multiplying the pencil~(\ref{eq3.13}) by the 
$
\bar M(x)=\mbox{diag} \{E_{d},M^{-1}(x),E_{\hat l}\}
$  
on the left, gives
\begin{equation}
\mbox{diag}\{E_{d}, \hat M(x),  E_{\hat l}\}+\lambda\ \mbox{diag}\{{\cal J}(x), E_{l},  {\cal N}(x)\},
\label{eq3.14}
\end{equation} 
with $\hat M(x)\equiv M^{-1}(x)- c E_{l}$. Consider the block $\hat M(x)$.  Its  characteristic equation takes the form 
$$
\det \left ( M(x)-\frac{1}{c+\zeta}E_{l}\right )=0,\ \  \zeta\neq -c\ \forall\ x\in {\cal U},
$$
where $\zeta\equiv\zeta(x)$ is some unknown function. From the latter it follows, that all eigenvalues of $\hat M(x)$ are identically equal to zero, because all eigenvalues of $M(x)$ are equal to $1/c$. Furthermore, the  rank of $A(x)$ is constant in ${\cal U}$. Taking into account the property of the rank, we conclude that the rank of  $\hat M(x)$ is independent of  $x\in {\cal U}$.  According to the lemma~\ref{l1}, there exists the nonsingular in ${\cal U}$ unitary matrix-function , satisfying the relation
$$
\bar U^{-1}(x) \hat M(x) \bar U(x)={\cal M}(x),
$$ 
where ${\cal M}(x)$ is some nilpotent matrix-function.

Let 
$
\hat U(x)\equiv\mbox{diag} \{E_{d}, \bar U(x), E_{\hat l}\}.
$
Multiplying the pencil~(\ref{eq3.14}) by the  matrix-functions $\hat U^{-1}(x)$ and $\hat U(x)$ on the left and on the right, respectively, we obtain the pencil~(\ref{eq3.1}).
Thus, we have proved the existence of the following  nonsingular matrix-functions in ${\cal U}$:
$$
P(x)=\hat U^{-1}(x)\bar M(x)\bar J(x)\tilde U^{-1}(x)T^{-1}(x)A_{1}^{-1}(x),
$$ 
$$
Q(x)=T(x)\tilde U (x)\hat U(x),
$$
which bring the pencil $A(x)+\lambda B(x)$ to the canonical form~(\ref{eq3.1}).  The elements of  $P(x)$ and $Q(x)$ belong to the spase ${C}^{p}({\cal U})$. Therefore the theorem is proved.
\end{pf}

In conclusion of this section, let us spend some lines to give a pair of remarks.
\begin{rmk}
Let us require, in circumstances of the theorem~\ref{t1}, the implementation of the following relations
\begin{equation}
{\rm rank}\ B(x)=\deg \left [\det(A(x)+\lambda B(x))\right ]
\label{eq3.15}
\end{equation}
and
\begin{equation}
{\rm rank}\ B(x)=\deg \left [\det(\mu A(x)+ B(x))\right ],
\label{eq3.16}
\end{equation}
where $\mu\equiv \mu(x)\in C^{p}({\cal U})$ is some unknown function, then the pencil $A(x)+\lambda B(x)$ is specified to be $p$-smoothly equivalent to the canonical form~(\ref{eq3.1}), in which
$$
{\cal N}(x)\equiv {\cal O}_{l}\ \mbox{and} \ {\cal M}(x)\equiv {\cal O}_{\hat l},
$$
where ${\cal O}_{l}$ is the zero block of order $l$. In this case $A(x)+\lambda B(x)$ is called the pencil satisfying the criterion ``rank-degree''. The structure of this pencil was investigated in~\cite{chist}.
\label{r1}
\end{rmk}

\begin{rmk}
If, in circumstances of the theorem~\ref{t1} and  conditions~(\ref{eq3.15}) and ~(\ref{eq3.16}), we additionally require that all the roots of the characteristic polynomial $\det(A(x)+\lambda B(x))$ are simple, then the pencil $A(x)+\lambda B(x)$  
is $p$-smoothly equivalent to a canonical form~(\ref{eq3.1}), in which
$$
{\cal N}(x)\equiv {\cal O}_{l},\ \ {\cal M}(x)\equiv {\cal O}_{\hat l},\ \ 
{\cal J}(x)={\rm diag} \left \{{-1}/{\lambda_{1}(x)},{-1}/{\lambda_{2}(x)},\dots,{-1}/{\lambda_{k}(x)}\right\},
$$
where $\lambda_{i}(x)$ are roots of the characteristic polynomial $\det(A(x)+\lambda B(x))$.
\label{r2}
\end{rmk}

\section{Examples}
\label{s4}
The goal of this section is to show the pair of simple examples to illustrate our theorem \ref{t1}.

{\sc Example 1.} Consider the pencil
\begin{equation}
A(x)+\lambda B(x)=
\left (
\begin{array}{ccc}
x_{1}+x_{2} & 0 & 0\\
0 & 0 & 0\\
0 & 0 & 1
\end{array} 
\right ) 
+\lambda 
\left (
\begin{array}{ccc}
1 & 0 & 0\\
0 & x_{1}x_{2} & -x_{2}^{2}\\
0 & x_{1}^{2} & -x_{1}x_{2},
\end{array} 
\right ) 
\label{eq4.1}
\end{equation} 
where 
$
x\in {\cal U}=[a, b]\times [a, b]\subseteq {\mathbb R}^{2},
$
and
$
 a>0.
$
Clearly, $\mbox{rank} \ A(x)=\mbox{rank} \ B(x)=2$  for $x\in {\cal U}$ and elements of $A(x)$ and $B(x)$ are  the analytic (polynomial) functions.
Characteristic polynomial for the pencil~(\ref{eq4.1}) is specified as
\begin{equation}
\det(A(x)+\lambda B(x))=x_{1}x_{2}\lambda (\lambda+x_{1}+x_{2}).
\label{eq4.2}
\end{equation} 
The roots of polynomial~(\ref{eq4.2})  are $\lambda_{1}\equiv 0$ and $\lambda_{2}(x)=-x_{1}-x_{2}$. They do not coincides in  ${\cal U}$ and their  multiplicity are constants in this domain. Furthermore, the leading coefficient of~(\ref{eq4.2}) is not equal to zero in ${\cal U}$. According to the theorem~\ref{t1}, the pencil~(\ref{eq4.1}) in ${\cal U}$ is 
equivalent to 
$$
\left (
\begin{array}{ccc}
1 & 0 & 0\\
0 & 0 & 0\\
0 & 0 & 1
\end{array} 
\right ) 
+\lambda 
\left (
\begin{array}{ccc}
{1}/{(x_{1}+x_{2})} & 0 & 0\\
0 & 1 & 0\\
0 & 0 & 0
\end{array} 
\right ). 
$$
The matrix-functions $P(x)$ and $Q(x)$ can be calculated step-by-step using proof of theorem~\ref{t1}. To do this, it is enough to put $c\equiv 1$.
We obtain these transforming matrices in the form
$$
P(x)=
\left (
\begin{array}{ccc}
1/(x_{1}+x_{2}) & 0 & 0\\
0 & 1/(x_{1}x_{2}) & 0\\
0 & -x_{1}/x_{2} & 1
\end{array} 
\right )
\ \mbox{and} \ 
Q(x)=
\left (
\begin{array}{ccc}
1 & 0 & 0\\
0 & 1 & x_{2}/x_{1}\\
0 & 0 & 1
\end{array} 
\right ). 
$$
One sees that the matrix-functions $P(x)$ and $Q(x)$ are nonsingular and analytic in ${\cal U}$.

{\sc Example 2.} Consider the pencil 
$$
A(x)+
\lambda\ B(x)=
\left (
\begin{array}{ccc}
0 & 0 & \gamma(x)\sigma (x)\\
0 & \gamma(x) & 0\\
0 & -\gamma(x)x_{2} & \gamma(x)
\end{array} 
\right )
$$
\begin{equation}
+\lambda 
\left (
\begin{array}{ccc}
0 & \sigma (x)x_{2} & -\upsilon (x)\sigma (x)\\
x_{1}^{2} & \upsilon (x) & -\upsilon (x)\sin(x_{1})\\
0 & -x_{1}(\upsilon (x)-1) & \upsilon (x)(\upsilon (x)-1)
\end{array} 
\right ), 
\label{eq4.3}
\end{equation} 
$$
x\in {\cal U}=[1, b]\times [1, b]\subseteq {\mathbb R}^{2},
$$
$$
\gamma(x)=x_{1}x_{2},\ \ \sigma (x)=x_{1}+x_{2},\ \ \upsilon (x)=x_{2}\sin(x_{1}).
$$
Observe, that 
there exists the minors of second order 
$$
A_{23}=
\left |\begin{array}{cc} 
\gamma(x) & 0\\ 
-\gamma(x) x_{2} & \gamma(x) 
\end{array}
\right |
\neq 0 \ \mbox{and} \ 
B_{23}=
\left |\begin{array}{cc} 
0 & \sigma (x)x_{2} \\ 
x_{1}^2 & \upsilon (x) 
\end{array}
\right |
\neq 0\ \ \forall\ x\in {\cal U}.
$$
Thus, $\mbox {rank} (A(x))=\mbox {rank} (B(x))=2$   $\forall\ x\in {\cal U}$. One sees that the characteristic polynomial for the pencil~(\ref{eq4.3}), that is 
$$
\det (A(x)+\lambda B(x))=-x_{1}^{4}x_{2}^{3}(x_{1}+x_{2})\lambda
$$
has single root $\lambda \equiv 0$ $\forall\ x \in {\cal U}$. Moreover, the leading coefficient of the characteristic polynomial is not equal zero in ${\cal U}$. According to our theorem~\ref{t1}, the pencil~(\ref{eq4.3}) in ${\cal U}$ is smoothly equivalent to
\begin{equation}
\left (
\begin{array}{ccc}
0 & 0 & 0\\
0 & 1 & 0\\
0 & 0 & 1
\end{array}
\right )
+\lambda
\left (
\begin{array}{ccc}
1 & 0 & 0\\
0 & 0 & \varphi (x)\\
0 & 0 & 0
\end{array}
\right ).
\label{eq4.4}
\end{equation}
The matrix-functions $P(x)$ and  $Q(x)$ are of the form
$$
P(x)=\frac{x_{1}^{2}}{\nu(x)x_{2}}
\left (
\begin{array}{ccc}
-\rho(x)\gamma(x) & \nu(x)x_{2} & \nu(x)\\
(\sin(x_{1})+x_{2})x_{1}^{2} & 0 & -\sin(x_{1})\sigma(x)x_{1}^{2}\\
(1-\upsilon(x))x_{1}^{2} & 0 & -\sigma(x)x_{1}^{2}
\end{array}
\right ),
$$
$$
Q(x)=\frac{1}{\rho(x)}
\left (
\begin{array}{ccc}
\rho(x) & 0 & 0\\ 
0 & \sin (x_{1}) & 1\\
0 & 1 & -\sin (x_{1})
\end{array}
\right )
$$
with $\nu(x)=\rho(x)\gamma(x)\sigma(x)$ and  $\rho(x)=(1+\sin^{2}(x_{1}))^{1/2}$. Multiplying the pencil~(\ref{eq4.3}) by the matrix-functions $P(x)$ and $Q(x)$ on the left and on the right, respectively, we obtain the function $\varphi(x)=(1+\sin^{2}(x_{1}))/x_{1}$. One sees that the elements of matrices  $P(x)$ and $Q(x)$ are analytic in ${\cal U}$.

In conclusion, let us observe, that the blocks ${\cal J}_{i}(x)$ in the canonical structure~(\ref{eq3.1}) are nonsingular matrix-functions of orders $p_{i}$, respectively. They do not have, generally speaking, Jordan structure. To make the blocks ${\cal J}_{i}(x)$ entirely coinciding with the Jordan blocks we must require additional conditions on the matrix-functions $A(x)$ and $B(x)$. In this paper we would not want to do this. Obtained the canonical form of the pencil~(\ref{eq3.1}) is quite sufficient to start the study of some class of the linear partial differential-algebraic equations and constructing numerical methods for them.





\bibliographystyle{model1a-num-names}


\begin{thebibliography}{99}

\bibitem{luch}
{\it W.~Lucht, K.~Strehmel and C.~Eichler-Liebenow} Indexes and special discretization methods for linear partial differential algebraic equations, BIT~39~(1999)~484--512.

\bibitem{deb}
{\it K.~Debrabant, K.~Strehmel} Convergence of Runge-Kutta methods applied to linear partial differential-algebraic equations, Applied Numerical Mathematics~53~(2005)~213--219.

\bibitem{god}
{\it C.\,K.~Godunov} Equations of the mathematical physics.  Nauka, Moscow, 1971 (in Russian).

\bibitem{gai1}
{\it S.\,V.~Gaidomak} Spline collocation method for linear singular hyperbolic systems, Comput. Math. Math. Phys.~48~(2008),~1161--1180.  

\bibitem{gai2}
{\it S.\,V.~Gaidomak} Three-layer finite-difference method for linear partial differential-algebraic systems,  
Differential Equations~46~(2010)~586--597.  

\bibitem{gai3}
{\it S.\,V.~Gaidomak} Stability of an implicit difference scheme for a linear differential-algebraic system of partial differential equations, 
Comput. Math. Math. Phys.~50~(2010)~673--683. 

\bibitem{chist}
{\it V.\,F.~Chistyakov} Algebraic-differential operators  with finite-dimensional kernel.  Nauka, Novosibirsk, 1996 (in Russian).

\bibitem{gant}
{\it F.\,R.~Gantmacher} The Theory of Matrices. Chelsea Publishing Company, New York 68, 1959. 

\bibitem{verb0}
{\it B.\,V.~Verbitskii} Some global property of matrix-functions that depend on several variables, Russian Math (Iz. Vuz)~22~(1978)~5--12.

\bibitem{verb1}
{\it B.\,V.~Verbitskii} A certain global property of matrix-valued functions that depend on several variables, Uspekhi Mat. Nauk~28~(1973)~233–-234 (in Russian). 

\bibitem{verb2}
{\it B.\,V.~Verbickii} On a global property of a matrix-valued function  of one variable, Math. USSR Sb.~20~(1973)~53--65. 

\bibitem{lank}
{P.~Lankaster} Theory of Matrices. Academic Press, New York-London, 1969. 

\bibitem{gin}
{\it H.~Gingold and P.\,F.~Hsieh} Globally analytic triangularization of a matrix function, Linear Algebra Appl.~169~(1992)~75–-101.

\bibitem{sil}
{\it L.\,M.~Silverman and R.\,S.~Bucy} Generalizations of a theorem of Dolezal, Theory of Computing Systems~4~(1969),~334--339.  

\bibitem{sib}
{\it P.\,F.~Hsieh and Y.~Sibuya} A global analysis of matrices of functions of several variables, J. Math. Anal. Appl.~14~(1966)~332–-340. 

\end{thebibliography}






\end{document}